\documentclass[11pt]{article}
\usepackage{amsmath}
\def\~{\tilde}
\def\D{{\Bbb D}}

\def\F{{\Bbb F}}
\def\Z{{\Bbb Z}}

\def\B{{\cal B}}

\def\N{{\cal N}}

\def\P{{\cal P}}

\def\Z{{\cal{Z}}}

\def\J{{\cal J}}

\def\ad{{\rm ad}}

\def\Proof{\noindent{\sl Proof.}\ }
\def\qed{{\hfill $\Box$ \medbreak}}
\usepackage{amssymb}

\newtheorem{defi}{Definition}[section]
\newtheorem{thm}[defi]{Theorem}
\newtheorem{lem}[defi]{Lemma}
\newtheorem{cor}[defi]{Corollary}

\newtheorem{prob}[defi]{Problem}

\newtheorem{prop}[defi]{Proposition}

\begin{document}

\title{Group identities on the units of algebraic algebras with applications to  restricted enveloping algebras}

\author{\sc Eric Jespers\thanks{Research
partially supported by Onderzoeksraad of Vrije
Universiteit, Fonds voor Wetenschappelijk Onderzoek
(Vlaanderen) and Flemish-Polish bilateral agreement
BIL2005/VUB/06.}, David Riley\thanks{Research
supported by NSERC of Canada.},  and Salvatore
Siciliano\thanks{The third author is grateful to the
Department of Mathematics at the University of
Western Ontario for its hospitality and support
while this research was carried out.}}

\date{}
\maketitle

\begin{abstract}
An algebra $A$ is called a GI-algebra if its group
of units $A^\times$ satisfies a group identity. We
provide positive support for the following two open
problems.
\begin{enumerate}
\item
Does every algebraic GI-algebra satisfy a polynomial
identity?
\item
Is every algebraically generated GI-algebra locally
finite?
\end{enumerate}
 \medbreak

 \noindent 2000 {\em Mathematics Subject Classification}.
 16R40; 16R50; 16U60; 17B35.
\end{abstract}

\section{Introduction}
Let $A$ be a unital associative algebra over a field
$\F$. We shall denote by $A^\times$ the group of all
multiplicative units of $A$. Recall that a group $G$
is said to satisfy a group identity whenever there
exists a nontrivial word $w(x_1,\ldots,x_n)$ in the
free group generated by $\{x_1,x_2,\ldots \}$ such
that $w(g_1,\ldots,g_n)=1$, for all $g_1,\ldots,g_n\in
A^\times$. By way of analogy with the custom of
referring to an algebra as a PI-algebra whenever it
satisfies a polynomial identity, an algebra $A$ is
sometimes called a GI-algebra if $A^\times$
satisfies a group identity. We shall also call a
Lie algebra satisfying a (Lie) polynomial identity
a PI-algebra.

The class of all GI-algebras has received
considerable attention recently. See \cite{DJJ} for
a recent and comprehensive overview. In particular,
it was shown in \cite{Liu00} that every algebraic
GI-algebra is locally finite. Furthermore, if the
base field is infinite then such an
algebra satisfies a polynomial identity (see also \cite{DG}).

In fact, rather
more can be said. But first we need to introduce
some more notation.

Recall that an associative algebra $A$ can be regarded
as a Lie algebra via
$[a,b]=ab-ba$, for every $a,b\in A$. As is customary,
we shall call a
Lie algebra bounded Engel if it satisfies a polynomial
identity of the form $$[x,y,\ldots,y]=0.$$ For
$x,y\in A^\times$, we write $(x,y)$ for the group
commutator $x^{-1}y^{-1}xy$. By analogy, a group
is called bounded Engel if it satisfies a group
identity of the form $$(x,y,\ldots,y)=1.$$

We shall
say that $A$ is Lie solvable (respectively,
bounded Engel or Lie nilpotent) to mean that $A$ is
solvable as a Lie algebra (respectively, bounded Engel or
nilpotent).

A polynomial identity is called non-matrix
if it is not satisfied by the algebra $M_2(\F)$ of
all $2\times 2$ matrices over $\F$.

For an algebra $A$, we denote by $\B(A)$ its prime
radical,  by $\J(A)$ its Jacobson radical, and by
$\Z(A)$ its center. Also, the set of all nilpotent
elements of $A$ will be denoted by $\N(A)$.

\begin{prop}
Let $A$ be an algebraic algebra over an infinite
field $\F$ of characteristic $p\geq 0$. Then the
following conditions are equivalent:
\begin{enumerate}
\item The algebra $A$ is a GI-algebra.
\item
The group of units $A^\times$ is solvable, in the
case when $p=0$, while $A^\times$ satisfies a group
identity of the form $(x,y)^{p^t}=1$ for some
natural number $t$, in the case when $p>0$.
\item
The algebra $A$ satisfies a non-matrix polynomial
identity.
\item
The algebra $A$ is Lie solvable, in the case when
$p=0$, while $A$ satisfies a polynomial identity of
the form $([x,y]z)^{p^t}=0$ for some natural number
$t$, in the case when $p>0$.
\end{enumerate}
Furthermore, in this case, $\N(A)=\B(A)$ is a
locally nilpotent ideal of $A$ and $A/\B(A)$ is both
commutative and reduced.
\end{prop}

\Proof
We need only collect various known results.
First we remark that every algebraic
algebra (or algebra generated by its algebraic elements)
over a field with
more than 2 elements is generated by its units. This
is a simple consequence of Wedderburn's theorem.
Consequently, if $A^\times$ satisfies a group identity
then Proposition 1.2 and Theorem 1.3 in
\cite{DJJ} apply yielding the fact that
$\N(A)$ is a locally nilpotent ideal of
$A$ coinciding with the prime radical $\B(A)$.
By Theorem 1.2 of \cite{BRT}, the subalgebra
$\F\cdot 1+\N(A)$ satisfies a non-matrix polynomial
identity. Since $A/\B(A)$ is commutative, by
Corollary 1.4 of \cite{DJJ}, it is easy to check
that $A$ itself satisfies a non-matrix polynomial
identity (see \cite{RW2}). The remaining
implications follow as in Theorems 1.3 and 1.4 in
\cite{BRT}.
\qed

In the previous proof we refer to \cite{DJJ} where algebras generated
by units are considered. Hence one might maybe expect that Proposition 1.1 is
valid in this context as well. However, this is not the case.
Indeed, let $\F G$ be the group algebra of a relatively free
nilpotent group $G$ of class $n>1$ over a field $\F$.
By  \cite{Seh} (Chapter V, Corollary 1.7),
all units of $\F G$ are of the form $\alpha g$, with $\alpha \in \F ^\times$ and $g\in G$.
In particular, the group of units of $\F G$ is nilpotent.
 On the other hand, in view of \cite{IP} and  \cite{Pass},
$\F G$ cannot satisfy
any polynomial identity.

Moreover, it is evident that Proposition 1.1
does not hold for finite base fields; however, the
following fundamental problem remains open.

\begin{prob}
Is every algebraic GI-algebra (over a finite
field) a PI-algebra?
\end{prob}

Problem 1.2 is a broad generalization of a
problem first posed by Brian Hartley. Hartley asked
if the group
algebra $\F G$ of a locally finite-$p$ group $G$
over a field $\F$ of characteristic $p>0$ is a
PI-algebra whenever it is a GI-algebra. This was
subsequently proved true for all periodic groups $G$
and fields $\F$ in \cite{GJV}, \cite{GSV},
\cite{Liu99} and \cite{LP}.

Understanding nilpotent algebras $N$ whose adjoint
group $1+N$ satisfies a group identity is crucial to
solving Problem 1.2 in general. Specifically, if
such an algebra $N$ satisfies a polynomial identity
that depends only on the particular group identity,
then Problem 1.2 would have a
general positive solution. Indeed, suppose that $A$ is any
algebraic algebra over a finite field $\F$ such that $A^\times$ satisfies
a particular group identity $w=1$. Let $H$ be an arbitrary finitely generated
(unital) subalgebra of $A$. Then, as mentioned earlier, $H$ is
finite-dimensional over $\F$, and hence actually finite.
Thus, by Wedderburn's theorem,
$H/\J(H)$ is the direct sum of matrix algebras (over finite fields).
So $(H/\J(H))^\times$
is the direct product of general linear groups each satisfying $w=1$.
But the degree of a general linear group
satisfying $w=1$ is bounded (see \cite{Liu00}). Hence, $H/\J(H)$
satisfies a polynomial identity determined only by $w=1$.
Thus, by our assumption applied to $N=\J(H)$,
$H$ itself would satisfy a polynomial identity
depending only on the group identity $w=1$.
Consequently, $A$ would be a PI-algebra.

Our present focus, however, is to give a positive
solution to Problem 1.2 for all algebraic algebras whose
units are either solvable or bounded Engel. Furthermore,
we intend to characterize these group-theoretic identities in
terms of the corresponding Lie polynomial identities.

We shall also address a second fundamental problem. Recall that an
 algebra is said to be algebraically generated
if it is generated by its algebraic elements.

\begin{prob}
Is every algebraically generated GI-algebra locally finite?
\end{prob}

Because algebraic GI-algebras are locally finite,
in order to solve Problem 1.3,  it suffices to show that every
algebraically generated GI-algebra is algebraic.
Even
in the case when the ground field is infinite, this
is not known. As evidence for a positive solution,
we offer the fact proved in \cite{DJJ} that every
group algebra of a periodically generated group over
an infinite field is locally finite provided its
group of units satisfies a group identity. We shall
offer presently some additional support to Problem 1.3 in the case
of nilpotently generated algebras. The case when
the base field is infinite was settled positively
by Theorem 1.2 in \cite{BRT}.

Finally, in the last section, we apply our general
results to the special case when the algebra in
question is the restricted enveloping algebra of a
restricted Lie algebra. Unlike the group algebra
case, very little is known about the group of units
of a restricted enveloping algebra.

\section{Algebraic GI-algebras over finite fields}

 We have seen in the Introduction that Proposition 1.1 answers Problem 1.2
in the case when the ground field is infinite. It is
clear, however, that Proposition 1.1 does not extend verbatim to
finite ground fields. We are unable to give a
complete solution to Problem 1.2 for finite fields,
but we will prove the case when the group of units
is either solvable or bounded Engel (see Corollary 2.3).
In fact,  we intend to give complete characterizations of these
conditions (for all but the smallest ground fields)
in Theorems 2.1 and 2.2.

\begin{thm}
Let $A$ be an algebraic algebra over a field $\F$
with $\vert \F \vert \ge 4$. Then $A^\times$ is
solvable if and only if all of the following conditions hold:
$A$ is Lie solvable,
$A/\J(A)$ is commutative, and there exists a chain
$0= J_0\subseteq J_1\subseteq \cdots \subseteq
J_m=\J(A)$ of ideals of $\J(A)$ such that every
factor $J_i/J_{i-1}$ is the sum of commutative
ideals of $J_m/J_{i-1}$.
\end{thm}
\medbreak \Proof First we prove necessity. Recall
from above that all algebraic GI-algebras are
locally finite. So, if $H$ is a finitely generated
(unital) subalgebra of $A$ and $\F$ is finite then
$H$ is actually finite. By Wedderburn's Theorem, it
follows that $H/\J(H)$ is a direct sum of matrix
algebras over field extensions of $\F$. Since $\vert
\F \vert\ge4$, the general linear group GL$_2(\F)$
is not solvable (see \cite{H}, Kapitel II, Satz 6.10).
It follows that $H/\J(H)$ is commutative and therefore
that $[H,H]H$ is nilpotent. This proves $[A,A]A$ is
locally nilpotent. In particular, $[A,A]A\subseteq
\J(A)$, so that $A/\J(A)$ is commutative. Because the
adjoint group
$1+\J(A)$ is solvable, by \cite{ASsolvable}, $\J(A)$
is Lie solvable and there exists a chain of ideals
of $\J(A)$ with the required properties. In
particular, it follows that $A$ itself is Lie
solvable.

It remains to consider the case when $\F$ is
infinite. Then, according to Proposition 1.1,
$\B(A)=\N(A)$ is locally nilpotent and $A/\B(A)$ is
commutative. In particular, $[A,A]A$ is locally
nilpotent and the proof follows as in the finite
field case.

Sufficiency follows easily from \cite{ASsolvable}.
\qed

Notice that the assumption on the minimal cardinality of
the ground field is required. Indeed, if
$A=M_2(\F_3)$, where $\F_3$ is the field of 3
elements, then $A^\times$=GL$_2(\F_3)$ is solvable
even though $A$ is not Lie solvable.

We now review some well-known facts for later use in the
proof of the next theorem.
Let $A$ be a finite-dimensional algebra over a field $\F$.
An element $x$ of $A$ is said to be semisimple if the
minimal polynomial of $x$ has no multiple roots
in any field extension of $\F$. In this case,
the adjoint map of $x$ turns out to be a semisimple
linear transformation of $A$. In particular, it is clear
that semisimple elements are central in a Lie nilpotent algebra.
Moreover, if $\F$ is perfect then the Wedderburn-Malcev Theorem implies
that for every element $x$ of $A$ there exist
$x_s,x_n\in A$ with  $x_s$ is semisimple and $x_n$ is
nilpotent such that $x=x_s+x_n$ and $[x_s,x_n]=0$.

We shall also use Lemma 2 from \cite{carter}. It can be stated as follows:
Let $A$ be a finite-dimensional algebra such that
$A/\J(A)$ is commutative and separable. If $G$ is a
nilpotent subgroup of $A^\times$ then the subalgebra
generated by $G$ is Lie nilpotent.

\begin{thm}
Let $A$ be an algebraic algebra over a perfect field
$\F\neq \F_2$.
\begin{enumerate}
\item
The group of units $A^\times$ is bounded Engel if and
only if $A$ is bounded Engel. In this case, $\N(A)$
is a locally nilpotent ideal such that
$A=\Z(A)+\N(A)$.
\item
The group of units $A^\times$ is nilpotent if and
only if $A$ is Lie nilpotent. Furthermore, the
corresponding nilpotency classes coincide.
\end{enumerate}
\end{thm}
\medbreak \Proof
 (1) Sufficiency holds for all rings (see \cite{RW1}).
For necessity, first recall that $A$ is locally
finite. Thus, in the case when $\F$ is finite,
if $H$ is a finitely generated (unital) subalgebra
of $A$ then $H$ is finite.
Thus, by Zorn's Theorem (see Theorem
12.3.4 in \cite{Rob}), the finite Engel group
$H^\times$ is nilpotent.
Another application of Wedderburn's Theorem easily
yields that $H/\J(H)$ is commutative. Since $\vert
\F\vert \ge3$, $H$ is generated by units. Also,
since $\F$ is perfect, $H/\J(H)$ is a separable
algebra. Therefore, Lemma 2 in \cite{carter} applies to $H$ and
hence $H$ is Lie nilpotent. Next we claim  that
$\J(H)=\N(H)$; that is, $\N(H)$ is an ideal of $H$.
Let $x,y\in \N(H)$. Since $H/\J(H)$ is commutative,
for every positive integer $t$, we have
$(x+y)^{p^t}\equiv x^{p^t}+y^{p^t}({\rm
mod}\,\J(H))$, so that, by the nilpotency of $x$ and
$y$ and the fact that $\J(H)$ is nilpotent, $x+y\in
\N(H)$. Also, if $a\in H$, from  $(ax)^{p^t}\equiv
a^{p^t}x^{p^t}({\rm mod}\,\J(H))$, it follows that
$ax\in\N(H)$, proving the claim. We conclude that
the set $\N(A)$ is a locally nilpotent ideal of $A$.
Next we prove the assertion that $A=\Z(A)+\N(A)$.
Indeed, let $x\in A$. Since $x$ is algebraic and
$\F$ is finite (and therefore perfect), by the remarks preceding the theorem
(applied to the associative subalgebra generated by $x$)
we have $x=x_s+x_n$, where $x_s$ is semisimple, $x_n$ is nilpotent, and
$[x_s,x_n]=0$.  Now let
$y\in A$ be arbitrary and set $B$ to be the
subalgebra generated by $x_s$ and $y$. Since (as
seen above) $B$ is a finite Lie nilpotent algebra,  the
semisimple elements of $B$ are central.  In particular,
$x_s$ and $y$ commute; hence, $x_s\in \Z(A)$ and
$A=\Z(A)+\N(A)$, as required. Finally, since the
adjoint group $1+\N(A)$ is bounded Engel, the locally nilpotent algebra
$\N(A)$ is a bounded Engel Lie algebra
(by \cite{ASengel}). This proves necessity in the
case when $\F$ is finite.

In the case when $\F$ is infinite, Proposition 1.1 informs us that
$\N(A)$ is a locally nilpotent ideal of $A$ such that $A/\N (A)$
is both commutative and reduced. Let $H$
be an arbitrary finitely generated (and therefore
finite-dimensional) subalgebra of $A$. It follows that $H/\N(H)$
is commutative, and thus $[H,H]H$ is nilpotent.
This implies that $H^\times$ is
solvable and so, by a well known result of Gruenberg
(\cite{G}), we see that the bounded Engel group
$H^\times$ is locally nilpotent. Now, since $\F$ is
perfect (by hypothesis) and $H$ is generated by its
group of units, Lemma 2 in \cite{carter} allows us to conclude that
$H$ is Lie nilpotent.  Thus, proceeding as above, we
obtain $A=\Z(A)+\N(A)$. Applying
\cite{ASengel} once more finishes the proof.

(2) Sufficiency holds for all rings (see \cite{GL}).
To prove necessity, we first infer from part (1)
that $\N(A)$ is a locally nilpotent ideal such that
$A=\Z(A)+\N(A)$.  But, according to a theorem of Du
(\cite{Du}), $\N(A)$ is Lie nilpotent of class
coinciding precisely with the nilpotency class of
the adjoint group $1+\N(A)$. The remaining
assertions follow at once. \qed

The ground field $\F_2$ was correctly omitted in
Theorem 2.2. Indeed, consider the restricted
enveloping algebra $u(L)$ of the restricted Lie
algebra $L$ over $\F_2$ with a basis $\{x,y\}$ such
that $[x,y]=x$, $x^{[2]}=0$, and $y^{[2]}=y$. Then
$u(L)^\times$ is isomorphic to the Klein four group.
Thus, $u(L)^\times$ is Abelian even though $u(L)$ is
not bounded Engel.

One can easily modify the proofs of Theorems 2.1 and 2.2
using the method outlined in the Introduction
in order to obtain the following result.

\begin{cor} Let $A$ be an algebraic algebra over
an arbitrary field $\F$. If $A^\times$ is either solvable or
bounded Engel  then $A$ is a PI-algebra.
\end{cor}

Now we turn to Problem 1.3.   An affirmative solution
(first proved in \cite{BRT}) follows
from Proposition 1.1 in the case when $A$ is generated by
nilpotent elements over an infinite field. A complete
solution will likely be difficult. However, for
GI-algebras generated by nilpotent elements over a
finite field, we offer the following result. Its proof uses
Theorem A from \cite{R2}, which
is essentially a corollary of a deep theorem due to
Zelmanov. Theorem A can be stated as follows:
Let $A$ be an associative algebra generated
 by a finite subset $X$. Suppose that the Lie subalgebra $L$ of $A$
generated by $X$ is a PI-algebra and
every Lie commutator (of length one or more)
in $X$ is nilpotent. Then $A$ is a nilpotent algebra.

\begin{thm}
Let $A$ be an algebra over a (finite) field of
characteristic $p>0$. Suppose $L\subseteq A$ is a
Lie subalgebra consisting of nilpotent elements. If
either $A^\times$ is solvable or bounded Engel then the
associative subalgebra $S$ generated by $L$ is
locally nilpotent.
\end{thm}
\Proof We can assume $L$ is finitely generated.
By Corollary 2 of \cite{BKL}, there
exists a positive integer $d$ such that $x^{d}\in
\B(A)$ for every $x\in \N(A)$. Clearly we can assume
$d=p^t$ for some $t$. Our
hypothesis now implies that every element in
the Lie subalgebra $L+\B(A)/\B(A)$ of $A/\B(A)$
is nilpotent of index at most $p^t$.
In particular, $L+\B(A)/\B(A)$ is bounded Engel (since $(\ad\, x)^{p^t}=\ad
(x^{p^t})$,
for each $x\in L$)
and thus $L+\B(A)/\B(A)$ satisfies a polynomial identity.
But $\B(A)\subseteq\J(A)$ is either Lie solvable or bounded Engel
(by \cite{ASsolvable} and \cite{ASengel}, respectively), and so
$L$ itself satisfies a polynomial identity. Consequently, $S$ is
nilpotent by Theorem A described above. \qed

\section{Applications to restricted enveloping algebras}
Let $L$ be a restricted Lie algebra over a field of
characteristic $p>0$. We recall that an element $x$
of $L$ is said to be $p$-nilpotent if there exists a
positive integer $n$ such that $x^{[p]^n}=0$. The
set of all $p$-nilpotent elements of $L$ will be
denoted by $\P(L)$. In particular, $L$ is called
$p$-nil if $L=\P(L)$. We say that $L$ is
algebraically generated if it is generated as a
restricted subalgebra by elements which are
algebraic with respect to the $p$-map of $L$. It
follows from the PBW theorem (see \cite{SF}, for
example) that the minimal polynomial of an element
$x$ in $L$ coincides precisely with the minimal
polynomial of $x$ when viewed as an element in the
restricted enveloping algebra $u(L)$ of $L$.

The main result of this section is the following. It
lends more support for the existence of an
affirmative solution to Problem 1.2.

\begin{thm}
Let $L$ be a restricted Lie algebra over an infinite
perfect field of characteristic $p>0$. If $L$ is
algebraically generated and $u(L)$ is a GI-algebra
then $u(L)$ is locally finite.
\end{thm}

\medbreak We split the proof into several lemmas.

\begin{lem}
Suppose that $L$ is a finite-dimensional restricted
Lie algebra over a perfect field $\F$ such that
$[L,L]\subseteq\P(L)$. Then the ideal $\P(L) u(L)$
coincides precisely with the set $\N(u(L))$. In
particular, $\B(u(L))=\P(L) u(L)$.
\end{lem}
\Proof Notice that since $\P(L)$ is a
finite-dimensional $p$-nil restricted ideal of $L$,
it follows that $\P(L) u(L)$ is a nilpotent ideal of
$u(L)$ (see \cite{RS}, for example).
As $u(L/\P(L))\cong u(L)/\P(L)u(L)$,  to complete
the proof it suffices to show that $u(L)$ is reduced
whenever $\P(L)=0$; consequently, we may also assume
that $L$ is Abelian. Suppose now,  to the contrary,
that $y$ is a nonzero nilpotent element in $u(L)$.
Fix an ordered basis $\{x_1,\ldots,x_n\}$ of $L$ and
write $y$ in its standard PBW representation: $$
y=\sum\alpha x_1^{a_1}\cdots x_n^{a_n}. $$ Then, by
the commutativity of $u(L)$, there exists a positive
integer $t$ such that $$ \sum\alpha^{p^t}
(x_1^{p^t})^{a_1}\cdots (x_n^{p^t})^{a_n}=y^{p^t}=0.
$$ It follows from the PBW theorem that
$\{x_1^{p^t},\ldots, x_n^{p^t}\}$ is a linearly
dependent set. Therefore, since $\F$ is perfect,
there exist $\beta_1,\ldots,\beta_n$ in $\F$, not
all zero, such that $$ \Big(\sum_{i=1}^n\beta_i
x_i\Big)^{p^t}=\sum_{i=1}^n\beta_i^{p^t}x_i^{p^t}=0.
$$ This contradicts our assumption that $L$ has no
nonzero $p$-nilpotent elements. \qed

Lemma 3.2 does not extend to arbitrary ground
fields. Indeed, let $\F$ be a field with positive
characteristic $p>0$ containing an element $\alpha$
with no $p$-th root in $\F$ and consider the Abelian
restricted Lie algebra $L=\F x+ \F y$ with
$x^{[p]}=x$ and $y^{[p]}=\alpha x$. Then it is easy
to check that $\P(L)=0$ while $0\neq
x^{p-1}y-y\in\N(u(L))$.

\begin{lem}
Suppose that $L$ is a finite-dimensional restricted
Lie algebra over an infinite field $\F$. If $u(L)$
is a GI-algebra then $[L,L]\subseteq \P(L)$.
\end{lem}
\Proof First note that since $L$ is
finite-dimensional, so is $u(L)$. Thus, by Proposition 1.1, $u(L)/\J(u(L))$ is
commutative. Since $\J(u(L))$ is nilpotent,
 it follows that every element in $[L,L]$ is $p$-nilpotent.
\qed

\begin{lem}
Let $L$ be an algebraically generated restricted Lie
algebra over an infinite perfect field $\F$ of
characteristic $p>0$ and suppose that $u(L)$ is a
GI-algebra. Then $\B(u(L))=\N(u(L))$ is locally
nilpotent. Furthermore,  $\P(L)=L\cap\B(u(L))$ is a
restricted ideal and $\B(u(L))=\P(L)u(L)$.
Consequently, $u(L)/\B(u(L)) \cong u(L/\P(L))$ is
reduced.
\end{lem}
\Proof As remarked in the proof of Proposition 1.1, such an
algebra $u(L)$ is generated by units. Hence, by
Theorem 1.3 in \cite{DJJ}, $\B(u(L))=\N(u(L))$ is a locally
nilpotent ideal of $u(L)$. It follows that
$\P(L)=L\cap\B(u(L))$ is a restricted ideal in $L$
such that $\P(L)u(L)\subseteq \B(u(L))$. In order to
prove the reverse inclusion, it suffices to assume
$\P(L)=0$. Let $\D(L)$ be the sum of all
finite-dimensional restricted ideals in $L$. Then,
according to Corollary 6.4 of \cite{BP}, $u(L)$ is
semiprime if and only if $u(\D(L))$ is
$L$-semiprime. But, under our assumption, Lemma 3.2
implies that $\B(u(I))=0$ for every
finite-dimensional restricted ideal $I$ of $L$.
Thus, $u(L)$ is semiprime, as required. \qed

\begin{lem}
Let $L$ be an algebraically generated restricted Lie
algebra over an infinite perfect field $\F$ of
characteristic $p>0$. If $u(L)$ is a GI-algebra then
$u(L)/\B(u(L))$ is commutative.
\end{lem}
\Proof We first apply the previous lemma in order to
pass to the case when $u(L)$ is reduced. Now let $x$
and $y$ be nonzero elements of $L$ with $x$
algebraic. We intend to show that necessarily
$[x,y]=0$.

Since $x$ is an algebraic element of $L$, there
exists a minimal positive integer $n$ with the
property that there exist $\alpha_0,\ldots,\alpha_n$
in $\F$, not all zero, such that
$$\sum_{i=0}^n\alpha_ix^{p^i}=0.$$ Clearly we may
also assume that $\alpha_n=1$. As explained above,
this is precisely the minimal polynomial of $x$ when
viewed as an algebraic element of $u(L)$. Since the
element $$[x,y]\sum_{i=0}^n\alpha_ix^{p^i-1}$$ has
square zero, it must be zero by assumption. Applying
the PBW theorem again yields that
$\{[x,y],x,x^p,\ldots,x^{p^{n-1}}\}$ is a linearly
dependent set in $L$. Consequently, the minimality
of $n$ implies $[x,y]\in \langle x\rangle_p$, where
$\langle x\rangle_p$ denotes the restricted
subalgebra generated by $x$. This implies
$[x,y,x]=0$ and hence $[y,x^p]=0$. Thus, we may
assume $\alpha _0=0$ for otherwise $x\in\langle
x^p\rangle_p$ commutes with $y$. In this case, the
element $$
 \sum_{i=1}^n\alpha_ix^{p^i-1}
$$ has square zero and hence, in fact, is zero
because $u(L)$ is reduced. However, this violates
our choice of $n$, and hence proves the lemma. \qed

Theorem 3.1 follows from a combination of Lemmas 3.4
and 3.5 and the fact that algebraic GI-algebras are
locally finite mentioned in the Section 1.

We conclude this section with the following corollaries, each of them is a
consequence of statements proved above and results
from \cite{RS} and \cite{RW2}.

\begin{cor}
Let $L$ be a restricted Lie algebra over a field
$\F$ of characteristic $p>0$. Then $u(L)$ is locally
finite  provided any one of the following conditions
hold.
\begin{enumerate}
\item The restricted Lie algebra $L$ is locally finite.
\item  The field $\F$ is both infinite and perfect, $L$ is
algebraically generated, and $u(L)$ is a GI-algebra.
\item
The field $\F$ is infinite, $L$ is generated by
$p$-nilpotent elements, and $u(L)$ is a GI-algebra.
\item
The field $\F$ is finite, $L$ is $p$-nil, and
$u(L)^\times$ is solvable or bounded Engel.
\end{enumerate}
\end{cor}

\begin{cor}
Let $L$ be a restricted Lie algebra over an infinite
field of characteristic $p>0$ and suppose that
$u(L)$ is algebraic. Then the following conditions
are equivalent.
\begin{enumerate}
\item The algebra $u(L)$ is a GI-algebra.
\item The algebra $u(L)$ satisfies a non-matrix polynomial identity.
\item
The restricted Lie algebra $L$ contains a restricted
ideal $I$ such that $L/I$ and $[I,I]$ are
finite-dimensional and $[L,L]$ is $p$-nil of bounded
index.
\end{enumerate}
\end{cor}

\begin{cor}
Let $L$ be a restricted Lie algebra over a field of
odd characteristic $p$ with at least 5 elements. If
$u(L)$ is algebraic then the following conditions
are equivalent.
\begin{enumerate}
\item The group of units $u(L)^\times$ is solvable.
\item The algebra $u(L)$ is Lie solvable.
\item
The derived subalgebra $[L,L]$ is both
finite-dimensional and $p$-nilpotent.
\end{enumerate}
\end{cor}

\begin{cor}
Let $L$ be a restricted Lie algebra over a perfect field
 of positive characteristic $p$ with at least 3
elements. If $u(L)$ is algebraic then the following
conditions are equivalent.
\begin{enumerate}
\item The group of units $u(L)^\times$ is bounded Engel.
\item The algebra $u(L)$ is bounded Engel.
\item
The restricted Lie algebra $L$ is nilpotent, $L$
contains a restricted ideal $I$ such that $L/I$ and
$[I,I]$ are finite-dimensional, and $[L,L]$ is
$p$-nil of bounded index.
\end{enumerate}
\end{cor}

\begin{cor}
Let $L$ be a restricted Lie algebra over a perfect field of
positive characteristic $p$ with at least 3
elements. If $u(L)$ is algebraic then the following
conditions are equivalent.
\begin{enumerate}
\item The group of units $u(L)^\times$ is nilpotent.
\item The algebra $u(L)$ is Lie nilpotent.
\item
The restricted Lie algebra $L$ is nilpotent and
$[L,L]$ is both finite-dimensional and
$p$-nilpotent.
\end{enumerate}
\end{cor}

Finally, we would like to mention that the Kurosh
problem for restricted enveloping algebras remains
open:

\begin{prob}
Is every algebraic restricted enveloping algebra
$u(L)$ locally finite?
\end{prob}
In other words, must the underlying restricted Lie
algebra $L$ be locally finite?

\bigbreak \noindent {ERIC JESPERS}\\ { Department of
Mathematics}\\ { Vrije Universiteit Brussel} \\ {
Pleinlaan 2, 1050 Brussel}\\ {Belgium}\\ {e--mail:
\it efjesper@vub.ac.be}

\bigbreak \noindent {DAVID RILEY}\\ { Department of
Mathematics}\\ { The University of Western Ontario}
\\ { London, Ontario, N6A 5B7}\\ {Canada}\\
{e--mail: \it dmriley\@@uwo.ca}

\bigbreak \noindent {SALVATORE SICILIANO}\\ {
Department of Mathematics ``E. De Giorgi''}\\ {
Universit\`{a} di Lecce} \\ { Via Prov. Lecce --
Arnesano, 73100, Lecce}\\ {Italy}\\ {e--mail: \it
salvatore.siciliano@unile.it}
\end{document}